\documentclass{amsart}
\usepackage{amssymb}
\def\version{Version W04-TsagasVolume-v2d.tex - Last Changed 13 November 2003}
\usepackage{a4}
\usepackage[notref,notcite]{showkeys}
\begin{document}
\date{\version}
\newtheorem{theorem}{Theorem}[section]
\newtheorem{example}{Example}
\font\pbglie=eufm10
\def\qedbox{\hbox{$\rlap{$\sqcap$}\sqcup$}}
\def\RR{{\text{\pbglie R}}}
\def\BB{{\mathcal{B}}}
\def\gg{{\text{\pbglie g}}}
\def\ffrac#1#2{{\textstyle\frac{#1}{#2}}}
\makeatletter
  \renewcommand{\theequation}{%
   \thesection.\alph{equation}}
  \@addtoreset{equation}{section}
 \makeatother

\title[Spectral geometry of the Riemann curvature operator]
{The spectral geometry of the Riemann curvature operator in the higher
signature setting}
\author{C. Dunn, P. Gilkey, R. Ivanova, and S. Nik\v cevi\'c}
\begin{address}{CD and PG: Mathematics Department, University of Oregon, Eugene Oregon 97403 USA.
  email:cdunn@darkwing.uoregon.edu and gilkey@darkwing.uoregon.edu}
\end{address}
\begin{address}{RI: Mathematics Department, University of Hawaii-Hilo, 200 West Kawili St.
Hilo, HI 96720-4091 USA. email: rivanova@hawaii.edu}\end{address}
\begin{address}{SN: Mathematical Institute, SANU, Knez Mihailova 35, p.p. 367, 11001 Belgrade, Yugoslavia. Email:
stanan@mi.sanu.ac.yu}
\end{address}
\begin{abstract} We study the spectral geometry of the Riemann curvature tensor
for pseudo-Riemannian manifolds and provide some examples illustrating the phenomena which can arise in the higher signature
setting.\\ {\bf Dedication:} This paper is dedicated to the memory of our colleague Prof. Gr.
Tsagas who studied the spectral geometry of the Laplacian.\end{abstract}
\maketitle
\section{introduction} Many authors have studied the spectral geometry of the Laplacian for a compact Riemannian manifold; see,
for example, Tsagas \cite{refTs94,refTs99} for additional references. Let
$\Delta_q:=(d\delta+\delta d)_q$ be the
$q$ form valued Laplacian. Tsagas \cite{TsTs91} showed that if $(M,g)$ is an $m$ dimensional closed Riemannian
manifold, then there exists $q=q(m)$ so that $(M,g)$ has the same $q$ spectrum as that of a round sphere $S^m$ if and only if
$(M,g)$ is in fact isometric to $S^m$. Thus round spheres are characterized by their $q$ spectrum; there are many other
results relating the spectrum of the Laplacian to the underlying geometry of the manifold.

In this brief note, we shall discuss the spectral geometry of the Riemann curvature tensor. Let $\nabla$ be the Levi-Civita
connection of a pseudo-Riemannian manifold $(M,g)$ of signature
$(p,q)$ and dimension $m:=p+q$. The associated curvature operator and curvature tensor are defined by setting:
\begin{eqnarray*}
&&R(x,y):=\nabla_x\nabla_y-\nabla_y\nabla_x-\nabla_{[x,y]},\quad\text{and}\\
&&R(x,y,z,w):=g(R(x,y)z,w)\,.
\end{eqnarray*}
Instead of examining the $L^2$ spectrum of the Laplacian, one
assumes that a natural operator which is associated to
$R$ has constant Jordan normal form on the associated domain of definition and then studies the attendant geometric
consequences. Such
questions were first raised in the seminal works of Ivanova-Stanilov
\cite{refRIGS} and Osserman \cite{refOss} and are a natural analogue of the questions studied for the
Laplacian. In the Riemannian setting ($p=0$), the operators in question are self-adjoint or
skew-adjoint and the spectrum determines the conjugacy class of the operator. However if the metric is
indefinite, then this is no longer the case, so one works with the Jordan normal form rather than with the
eigenvalues.

Here is a brief outline to this paper. In Section \ref{sect-2}, we present the basic definitions and give a short survey of the
current state of the field. Certain questions are essentially settled for Riemannian manifolds or for
Lorentzian manifolds ($p=1$). However, relatively little is known in the higher signature setting. In
Section
\ref{sect-3} we shall exhibit some examples to illustrate phenomena which arise for manifolds of higher signature. We
conclude the paper with a rather lengthy bibliography to serve as a partial introduction to the field.

\section{Natural operators defined by the Riemann curvature tensor}
\label{sect-2}

\subsection{The Jacobi operator $J$}  Set
\begin{equation}\label{eqn-2.a}
J(x):y\rightarrow R(y,x)x\,;
\end{equation}
this self-adjoint operator plays an important role in
the study of geodesic sprays. One says
$(M,g)$ is {\it spacelike} (resp. {\it timelike}) {\it Jordan Osserman} if the Jordan normal form of $J$ is constant on
the bundle of  unit spacelike (resp. timelike) tangent vectors in $TM$. 

In the Riemannian setting $(p=0)$, the Jordan normal form is determined by the eigenvalue
structure and, as every vector is spacelike, we shall drop the qualifiers
`spacelike' and `Jordan'. This is not true in the higher signature context
which is why we focus on the Jordan normal form, i.e. the conjugacy class, instead of only
on the eigenvalue structure.

One has the following results due to Chi \cite{refChi} and
to Nikolayevsky
\cite{Ni03,Ni04} in the Riemannian setting
 and to N. Bla\v zi\'c et. al. \cite{refBBG} and to Garc\'{\i}a--R\'{\i}o et. al. \cite{refGKV} in the Lorentzian setting; the
classification is essentially complete here:
\begin{theorem}\label{thm-2.1}\ \begin{enumerate}
\item Let $(M,g)$ be an Osserman Riemannian manifold of dimension $m\ne8,16$. Then either $(M,g)$ is locally
isometric to a rank $1$ symmetric space or $(M,g)$ is flat.
\item Let $(M,g)$ be a spacelike or timelike Jordan Osserman Lorentzian
manifold. Then $(M,g)$ has constant sectional curvature.
\end{enumerate}
\end{theorem}
We refer to
\cite{GKV02} for further details; it contains an excellent discussion of the spectral geometry of the
Jacobi operator.

\subsection{The higher order Jacobi operator} Let
$S(\pi)$ be the sphere of unit spacelike (resp. timelike) vectors in a spacelike (resp. timelike) $k$ dimensional subspace of
$TM$. Let
$$J(\pi):=\textstyle\ffrac1{\operatorname{vol}\{S(V)\}}\int_{x\in S(V)}J(x)dx$$
be the average Jacobi operator. If $\{e_1,...,e_k\}$ is an orthonormal basis for $\pi$, then modulo a suitable normalizing
constant which plays no role in our development,
$$J(\pi)=\textstyle\sum_iJ(e_i)\,.$$
These operators were first defined by Stanilov and Videv 
\cite{refSV} in the Riemannian context. One says that $(M,g)$ is {\it spacelike} (resp. {\it timelike}) {\it Jordan $k$-Osserman}
if the Jordan normal form of
$J$ is constant on the Grassmannian of spacelike (resp. timelike) $k$ planes of $TM$, where $2\le k\le q$ (resp. 
$2\le k\le p)$. One has the following classification result in the
Riemannian
\cite{refGi01} and Lorentzian
\cite{refGiSt02} settings:

\begin{theorem}\label{thm-2.2} Let $(M,g)$ be a spacelike Jordan $k$-Osserman pseudo-Riemannian manifold of signature $(p,q)$.
\begin{enumerate}
\item If $p=0$ and if $2\le k\le m-2$, then $(M,g)$ has constant sectional curvature.
\item If $p=0$ and if $k=m-1$, then either $(M,g)$ is locally isometric to a rank $1$ symmetric space or $(M,g)$ is flat.
\item If $p=1$ and if $2\le k\le m-1$, then $(M,g)$ has constant sectional curvature.
\end{enumerate}
\end{theorem}

\subsection{The skew-symmetric curvature operator $\mathcal{R}$} Let $\{e_1,e_2\}$ be an orthonormal basis for an oriented
spacelike or timelike $2$ plane $\pi$. One defines
\begin{equation}\label{eqn-2.b}
\mathcal{R}(\pi):=R(e_1,e_2);
\end{equation} 
this skew-symmetric operator is independent of the particular
orthonormal basis which was chosen for $\pi$. It is, however, sensitive to the orientation of $\pi$; if $-\pi$ denotes
$\pi$ with the opposite orientation, then $\mathcal{R}(-\pi)=-\mathcal{R}(\pi)$. The manifold $(M,g)$ is said to be {\it
spacelike} (resp. {\it timelike}) {\it Jordan Ivanov-Petrova} if for every $P\in M$, the Jordan normal form of $\mathcal{R}$ is
constant on the spacelike (resp. timelike) $2$ planes in
$T_PM$; in contrast to the Jacobi operator, the Jordan normal form is allowed to vary with the point in question.

In addition to manifolds of constant sectional curvature, there are warped product metrics which are both spacelike and
timelike Jordon Ivanov-Petrova. Let
$M:=I\times N$, where $I$ is an open sub-interval of $\mathbb{R}$ and where $(N,g_N)$ is a pseudo-Riemannian
manifold with constant sectional curvature
$\kappa$. Let
\begin{equation}\label{eqn-2.c}
g_M:=\varepsilon dt^2+\{\varepsilon\kappa t^2+At+B\}g_N\quad\text{for}\quad\varepsilon=\pm1\,.
\end{equation}
The sub-interval $I$ is chosen so that the warping function $\kappa t^2+At+B\ne0$ (this ensures the metric is non-degenerate) and
so that
$A^2-\varepsilon\kappa B\ne0$ (this ensures the metric is not flat). The manifold $(M,g_M)$ is then both spacelike and timelike
Jordan Ivanov-Petrova
\cite{refGil}.

 The classification of such manifolds is essentially complete in the Riemannian
setting \cite{refGi99,refGLS,IvPe} and in the Lorentzian setting \cite{GiZa,refZ}; there are some partial results due
to Stavrov
\cite{St03} if
$p\ge2$:

\begin{theorem}\label{thm-2.3} Let $(M,g)$ be a pseudo-Riemannian manifold of
signature $(p,q)$ which is spacelike Jordan Ivanov-Petrova.
\begin{enumerate}
\item If $p=0$, if $m\ge4$, and if $m\ne7$, then either $(M,g)$ has constant sectional curvature or $(M,g)$ is locally isometric
to a manifold as in Equation (\ref{eqn-2.c}).
\item If $q\ge11$ and if $\{q,q+1\}$
 does not contain a power of $2$, then either $(M,g)$ has constant sectional curvature or $(M,g)$ is locally isometric
to a manifold as in Equation (\ref{eqn-2.c}).
\item If $q\ge11$, if $1\le p\le(q-6)/4$, if the set $\{q,q+1,...,q+p\}$ does not contain a
power of $2$, and if $\mathcal{R}(\pi)$ is not nilpotent for some $2$ plane $\pi$, then either $(M,g)$ has constant sectional
curvature or $(M,g)$ is locally isometric to a manifold as  in Equation (\ref{eqn-2.c}).
\end{enumerate}
\end{theorem}

Assertion (3) can be used to derive results for timelike Jordan Ivanov-Petrova manifolds by changing the sign of the metric in
question and interchanging the roles of $p$ and $q$. We refer to \cite{refGil} for additional results concerning
Ivanov-Petrova manifolds. The manifolds given in Sections \ref{subs-a} and \ref{subs-3.4} will be spacelike Jordan
Ivanov-Petrova but will neither have constant sectional curvature nor have the form given in Equation
(\ref{eqn-2.c}).

\subsection{The Stanilov operator} Let $\pi$ be a spacelike or timelike $2$ plane in $TM$. Let $\operatorname{Gr}_2(\pi)$
be the Grassmannian of oriented
$2$ planes in
$\pi$. The Stanilov operator is an average of the square of the skew-symmetric curvature operator:
$$\Theta(\pi):=\ffrac1{\operatorname{vol}(Gr_2(\pi))}\textstyle\int_{\pi\in\operatorname{Gr}_2(\pi)}\mathcal{R}(\pi)^2d\pi\,.$$
This operator is self-adjoint; it is necessary to square $\mathcal{R}$ to obtain a non-zero average since
$\mathcal{R}(-\pi)=-\mathcal{R}(\pi)$.

If $\{e_1,...,e_k\}$ is an orthonormal basis for $\pi$, then modulo a suitable normalizing constant
which plays no role in the development,
$$\Theta(\pi)=\textstyle\sum_{i<j}\mathcal{R}(e_i,e_j)^2\,.$$
This operator was first defined by Stanilov \cite {Stx,Sty} in the Riemannian context.
One says that $(M,g)$ is {\it spacelike} (resp. {\it timelike}) {\it Jordan $k$-Stanilov} if the Jordan normal form of
$\Theta$ is constant on the Grassmannian of spacelike (resp. timelike) $k$ planes of $T_PM$, where $2\le k\le q$ 
(resp. $2\le k\le p)$ for every $P\in M$; as with the skew-symmetric curvature operator, the Jordan normal form is permitted to
vary with the point of
$M$.

We refer to \cite{refGNV} for the proof of the following result:
\begin{theorem}\label{thm-2.4}
Let $(M,g)$ be a connected spacelike Jordan
Ivanov-Petrova pseudo-Riemannian manifold of signature $(p,q)$. Assume either that $(p,q)=(0,4)$ or that $q\ge5$. Assume
that
$R(\pi)$ is not nilpotent for at least one spacelike $2$ plane in $TM$ and that $\mathcal{R}$ has spacelike rank $2$ for
all $P\in M$. Then
\begin{enumerate}
\item $(M,g)$ is spacelike Jordan $k$-Stanilov for any $2\le k\le q$.
\item $(M,g)$ is timelike Jordan $k$-Stanilov for any $2\le k\le p$.
\end{enumerate}
\end{theorem}

\begin{theorem}\label{thm-2.5} Let $(M,g)$ be a connected Riemannian manifold of dimension $m$, where $m\ne 3,7$.
If $(M,g)$ is $2$-Stanilov, then $(M,g)$ is Ivanov-Petrova.
\end{theorem}

We also refer to \cite{TzVi99} for related results. We shall exhibit manifolds which are spacelike $k$ Jordan
Stanilov for any $k$ but which are not spacelike Jordan Ivanov-Petrova in the higher signature setting in Section \ref{sect-3}.

\subsection{Subspaces of mixed type} It is not necessary to restrict to spacelike or timelike planes in discussing the
skew-symmetric curvature operator, the higher order Jacobi operator, or the Stanilov operator. Let
$\{e_1,e_2\}$ be a oriented basis for a non-degenerate oriented
$2$ plane
$\pi$ of signature
$(1,1)$. Let
$g_{ij}=g(e_i,e_j)$ describe the metric on $\pi$. Set:
$$\mathcal{R}(\pi):=|\det(g)|^{-1/2}R(e_i,e_j)\,.$$
One says that $(M,g)$ is {\it mixed Jordan Ivanov-Petrova} if the Jordan normal form of $\mathcal{R}(\pi)$ is constant on the
Grassmannian of oriented $2$ planes of signature $(1,1)$ in $T_PM$ for every $P\in M$; the Jordan normal form is allowed to vary
with the point. The manifolds described in Equation (\ref{eqn-2.c}) are spacelike, timelike, and mixed Jordan Ivanov-Petrova. In
Theorem \ref{thm-3.2x}, we present manifolds which are spacelike and timelike Jordan Ivanov-Petrova but
which are not mixed Jordan Ivanov-Petrova. In Theorem \ref{thm-3.6}, we present manifolds which are
spacelike Jordan Ivanov-Petrova but which are neither mixed nor timelike Jordan Ivanov-Petrova. Thus
these concepts are distinct.

Similarly, let $\{e_1,...,e_k\}$ be a basis for an unoriented non-degenerate
$k$ plane
$\pi$. Let
$g_{ij}:=g(e_i,e_j)$. One can define
\begin{eqnarray*}
&&J(\pi)(y):=\textstyle\sum_{ij}g^{ij}R(y,e_i)e_j,\quad\text{and}\\
&&\Theta(\pi)y:=\textstyle\sum_{ijkl}g^{ij}g^{kl}R(e_i,e_k)R(e_j,e_l)y\,.
\end{eqnarray*}
Using these operators, one can define the notions {\it Jordan Osserman of type $(r,s)$} and {\it Jordan Stanilov of type $(r,s)$}
by requiring that the Jordan normal form of $J(\pi)$ or $\Theta(\pi)$ is constant on the Grassmannian of non-degenerate planes of
type $(r,s)$. We omit details in the interests of brevity.

\subsection{Conformal geometry} Let $\{e_1,...,e_m\}$ be a basis for $T_PM$. Let
$$
\rho(x,y):=\textstyle\sum_{ij}g^{ij}R(x,e_i,e_j,y)\quad\text{and}\quad
\tau:=\textstyle\sum_{ij}g^{ij}\rho(e_i,e_j)
$$ 
be the {\it Ricci tensor} and the {\it scalar curvature}, respectively. The {\it Weyl conformal curvature} $W$ is
then defined by setting:
\begin{eqnarray*}
&&W(x,y,z,w):=R(x,y,z,w)+\ffrac1{(m-1)(m-2)}\tau\{g(x,w)g(y,z)-g(x,z)g(y,w)\\
&-&\ffrac1{m-2}\{\rho(x,w)g(y,z)+\rho(y,z)g(x,w)-\rho(x,z)g(y,w)-\rho(y,w)g(x,z)\}\,.
\end{eqnarray*}
One generalizes Equations (\ref{eqn-2.a}) and (\ref{eqn-2.b}) setting:
$$
J_W(x):y\rightarrow W(y,x)x\quad\text{and}\quad\mathcal{R}_W(\pi):y\rightarrow W(e_1,e_2)y\,.
$$
One says that $(M,g)$ is {\it conformally spacelike} (resp. {\it timelike}) {\it Jordan Osserman} if the Jordan
normal form of $J_W$ is constant on $S(T_PM)$; the Jordan normal form is permitted to vary with the point $P\in M$.
Similarly, one says that
$(M,g)$ is {\it conformally spacelike} (resp. {\it timelike}) {\it Ivanov-Petrova} if the Jordan normal form of
$\mathcal{R}_W$ is constant on the appropriate Grassmannian of $T_PM$; again, the Jordan normal form is permitted to
vary with the point of $M$. 

One says that two metrics $g_1$ and $g_2$ on $M$ are {\it conformally equivalent} if there exists a smooth positive conformal
factor
$\alpha\in C^\infty(M)$ so $g_1=\alpha g_2$.
One then has
$$W_{g_1}=\alpha W_{g_2}\,.$$
The notions defined above are conformal invariants \cite{refBGNS}:
\begin{theorem}\label{thm-2.6} Let $g_1$ and $g_2$ be conformally equivalent pseudo-Riemannian metrics on a manifold $M$. Then:
\begin{enumerate}
\item $(M,g_1)$ is conformally spacelike (resp. timelike)
Jordan Osserman if and only if $(M,g_2)$ is conformally spacelike (resp. timelike) Jordan Osserman.
\item $(M,g_1)$ is conformally spacelike (resp. timelike)
Jordan Ivanov-Petrova if and only if $(M,g_2)$ is conformally spacelike (resp. timelike) Jordan Ivanov-Petrova.
\end{enumerate}
\end{theorem}

One also has
\begin{theorem}\label{thm-2.7}
 If $(M,g)$ is Einstein, then $(M,g)$ is conformally spacelike (resp. timelike) Jordan Osserman if and only if
$(M,g)$ is pointwise spacelike (resp. timelike) Jordan Osserman.
\end{theorem}

The classification of conformally spacelike Jordan Osserman manifolds is complete in certain settings. One says
that $(M,g)$ is {\it conformally flat} if $W=0$ or, equivalently, if $g$ is locally conformally equivalent to a flat
metric. Note that metrics of constant sectional curvature are conformally flat. We refer to
\cite{refBGNS} for the proof of Assertion (1) and to
\cite{refBG03} for the proof of Assertion (2) in the following result:

\begin{theorem}\label{thm-2.8} Let $(M,g)$ be a conformally Osserman Riemannian manifold of dimension $m$.
\begin{enumerate}
\item If $m$ is odd, then $(M,g)$ is conformally flat.
\item If $m\equiv2$ mod $4$ and if $P$ is a point of $M$ where $W_P\ne0$, then $(M,g)$ is locally
conformally equivalent near $P$ either to complex projective space with the Fubini-Study metric or to the negative curvature dual.
\end{enumerate}
\end{theorem}

\begin{theorem}\label{thm-2.9} If $(M,g)$ is a conformally spacelike or conformally timelike Jordan Osserman Lorentzian manifold,
then $(M,g)$ has constant sectional curvature.
\end{theorem}

Similarly, there are results for conformally spacelike Jordan Ivanov-Petrova manifolds \cite{refBGNS}:
\begin{theorem}\label{thm-2.10}\ \begin{enumerate}
\item Let $(M,g)$ be a conformally Ivanov-Petrova Riemannian manifold of dimension $m\ne 3,7$. Then $(M,g)$ is
conformally flat.
\item Let $(M,g)$ be a connected pseudo-Riemannian manifold of signature $(p,q)$ which is conformally spacelike Ivanov-Petrova. Assume that
$q\ge11$, that
$p\le\frac{q-6}4$, and that $\{q,q+1,...,q+p\}$ does not contain a power of $2$. Then either $W(\pi)$ is nilpotent for every spacelike $2$
plane or $(M,g)$ is conformally flat.
\end{enumerate}
\end{theorem}

\section{The higher signature setting}\label{sect-3}

As noted in the previous section, there are many classification results available in the Riemannian and Lorentzian settings. The
situation is much less clear in the higher signature setting ($p>1,q>1$). 

\subsection{Curvature homogeneous manifolds} We follow Kowalski, Tricerri, and
Vanhecke \cite{KTV91,KTV92} and say that $(M,g)$ is {\it curvature homogeneous} if
given any two points
$P,Q\in M$, there is an isomorphism $\Psi:T_PM\rightarrow T_QM$ so $\Psi^*g_Q=g_P$
and so $\Psi^*R_Q=R_P$. Similarly, $(M,g)$ is said to be
{\it locally homogeneous} if the local isometries of $(M,g)$ act transitively on $(M,g)$. 

The manifold $(M,g)$ is said to be {\it locally
symmetric} if $\nabla R=0$. Locally symmetric manifolds are locally homogeneous and locally  homogeneous manifolds are curvature
homogeneous. What is interesting from our point of view is that the converse fails in general; there are curvature homogeneous
manifolds which are not locally homogeneous.

There is by now an extensive literature on the subject of curvature homogeneous manifolds in the Riemannian setting, see, for
example, the discussion in \cite{BV98,V91}. There are also a number of papers in the Lorentzian setting
\cite{BM00,BV97} and also in the affine setting \cite{O96}. There are, however, few papers in
the higher dimensional setting -- and those that exist appear in the study of $4$ dimensional neutral signature
Osserman manifolds, see, for example, \cite{BCR01}. In this section we present two families of examples
which illustrate many of the phenomena which can arise. In Section \ref{subs-a}, we exhibit pseudo-Riemannian manifolds in
balanced neutral signature \cite{refDG,GIZ02,GIZ03} and of signature $(2s,s)$ \cite{GN03}. We refer to
\cite{BBZ01,GVV98} for other examples and to
\cite{GKV02} for a more complete bibliography.

\subsection{Algebraic curvature tensors} It is convenient at this stage to introduce a purely algebraic formalism. Let $V$ be an
$m$ dimensional finite dimension real vector space. We say that $A\in\otimes^4V^*$ is an {\it algebraic curvature tensor} if $A$
satisfies the usual symmetries of the Riemann curvature tensor:
\begin{eqnarray*}
A(x,y,z,w)=-A(y,x,z,w)=A(z,w,x,y),\\
A(x,y,z,w)+A(y,z,x,w)+A(z,x,y,w)=0\,.
\end{eqnarray*}
Let $g_V$ be a non-degenerate symmetric inner product of signature $(p,q)$ on $V$. The associated curvature operator is then
characterized by the identity:
$$g_V(A(x,y)z,w)=A(x,y,z,w)\,.$$

Consider a triple $\mathcal{V}:=(V,g_V,A)$, where $A$ is an algebraic curvature tensor on $V$ and $g$ is an inner
product on $V$. We say that $\mathcal{V}$ is a {\it model} for a pseudo-Riemannian manifold $(M,g)$ if given any
point $P\in M$, there is an isomorphism $\phi_P:T_PM\rightarrow\mathcal{V}$ so that $\phi_P^*g_V=g|_{T_PM}$ and
$\phi_P^*A=R_P\in\otimes^4T_P^*M$. Clearly $(M,g)$ is curvature homogeneous if and only if there exists a model
$\mathcal{V}$ for $(M,g)$.

\subsection{Signature $(p,p)$}\label{subs-a} We follow the discussion in \cite{refDG,GIZ02,GIZ03}. Suppose $p\ge3$
henceforth.
 Introduce coordinates $(x,y)=(x_1,...,x_p,y_1,...,y_p)$ on
$\mathbb{R}^{2p}$. Let $\mathcal{O}$ be an open subset of $\mathbb{R}^p$ and let $f=f(x)\in C^\infty(\mathcal{O})$. We define a
non-degenerate pseudo-Riemannian metric $g_f$ of balanced signature
$(p,p)$ on
$M:=\mathcal{O}\times\mathbb{R}^p$ by setting:
$$
g_f(\partial_i^x,\partial_j^x)=\partial_i^xf\cdot\partial_j^xf,\quad
g_f(\partial_i^x,\partial_j^y)=g_f(\partial_j^y,\partial_i^x)=\delta_{ij},\quad\text{and}\quad
g_f(\partial_i^y,\partial_j^y)=0\,.
$$

This pseudo-Riemannian metric arises as a hyper-surface metric. Let 
$$\{\alpha_1,...,\alpha_p,\beta_1,...,\beta_p,\gamma\}$$
be a basis for $\mathbb{R}^{2p+1}$. Define an inner product of signature $(p,p+1)$ on $\mathbb{R}^{2p+1}$ whose non-zero
components are given, up to the usual $\mathbb{Z}_2$ symmetries, by
$$
\langle\alpha_i,\beta_j\rangle=\delta_{ij}\quad\text{and}\quad\langle\gamma,\gamma\rangle=1\,.
$$
We define an embedding of $\mathcal{O}\times\mathbb{R}^p$ as a graph in $\mathbb{R}^{2p+1}$ by
setting:
$$\Psi_f(x,y):=x_1\alpha_1+...+x_p\alpha_p+y_1\beta_1+...+y_p\beta_p+f(x)\gamma\,.$$
It is then immediate that $\Psi_f^*\langle\cdot,\cdot\rangle=g_f$. The normal to the surface is:
$$\nu(x,y)=-(\partial_1^xf)\beta_1-...-(\partial_p^xf)\beta_p+\gamma\,.$$
Let $H_{ij}=\partial_i^x\partial_j^xf\in M_p(\mathbb{R})$ be the Hessian. The second fundamental form $L$
and curvature tensor are given by:
\begin{eqnarray*}
&&L(\partial_i^x,\partial_j^x)=H_{ij},\quad,
L(\partial_i^x,\partial_j^y)=L(\partial_j^y,\partial_i^x)=L(\partial_i^y,\partial_j^y)=0,\\
&&R(x,y,z,w)=L(x,w)L(y,z)-L(x,z)L(y,w)\,.
\end{eqnarray*}
The only non-zero action of the curvature operator is:
$$
R(\partial_i^x,\partial_j^x):\partial_k^x\rightarrow R_{ijkl}\partial_l^y\,.
$$
The curvature operator is $2$-nilpotent as:
\begin{eqnarray*}
&&R(\cdot,\cdot):\operatorname{Span}\{\partial_i^x\}\rightarrow\operatorname{Span}\{\partial_i^y\},\\
&&R(\cdot,\cdot):\operatorname{Span}\{\partial_i^y\}\rightarrow\{0\}\,.
\end{eqnarray*}
As the Ricci operator is nilpotent, $(M,g_f)$ is Ricci flat and Einstein.

Let $\{X_1,...,X_p,Y_1,...,Y_p\}$ be a basis for $V:=\mathbb{R}^{2p}$. We define an inner product $g_V$ on $V$ of
signature $(p,p)$ and an algebraic curvature tensor whose non-zero components are, up to the usual $\mathbb{Z}_2$ symmetries:
\begin{eqnarray*}
&&g_V(X_i,Y_j)=g_V(Y_j,X_i)=\delta_{ij},\quad\text{and}\\
&&A_V(X_i,X_j,X_k,X_l)=\delta_{il}\delta_{jk}-\delta_{ik}\delta_{jl}\,.
\end{eqnarray*}
Let $\mathcal{V}_{p,p}:=(V,g_V,A_V)$. We then have \cite{refDG}
\begin{theorem}\label{thm-3.1} Let $p\ge2$. Assume that $H$ is positive definite on $\mathcal{O}$.
\begin{enumerate}
\item $(M,f)$ is curvature homogeneous with model $\mathcal{V}_{p,p}$. 
\item If $p\ge3$, then $(M,g_f)$ is not locally homogeneous
for generic $f$.
\end{enumerate}\end{theorem}

If $f(x)=x_1^2+...+x_p^2$, then $(M,g_f)$ is a local symmetric space; thus $\mathcal{V}$ is the model space of a symmetric space
of signature $(p,p)$. These manifolds form a nice family of examples to study the spectral geometry of the Riemann curvature
tensor. One has
\cite{GIZ03}:

\begin{theorem}\label{thm-3.2x} Assume that $H$ is non-degenerate. 
 \begin{enumerate}
\item The manifold $(M,g_f)$ is spacelike and timelike Jordan Ivanov-Petrova.
\item The manifold $(M,g_f)$ is not mixed Jordan Ivanov-Petrova.
\item If $p=2$ or if 
$p\ge3$ and if $H$ is definite, then $(M,g_f)$ is spacelike and timelike Jordan Osserman
\item If $p\ge3$ and if $H$ is indefinite, then $(M,g_f)$ is neither
 spacelike Jordan Osserman nor timelike Jordan Osserman.
\end{enumerate}\end{theorem}

Since $R(\cdot,\cdot)^2=0$, the manifolds $(M,g_f)$ are spacelike and timelike Jordan $k$-Stanilov for $2\le k\le p$. If $H$ is
degenerate, then $(M,g_f)$ is neither spacelike nor timelike Jordan Ivanov-Petrova. Thus there are examples of manifolds
which are Jordan $k$-Stanilov but not Jordan Ivanov-Petrova.

Let $\mathbb{R}^{(a,b)}$ be a flat manifold with a metric of signature $(a,b)$. 
One has \cite{GIZ03}:
\begin{theorem}\label{thm-3.3}  Assume $H_f$ is definite.
Let $N:=M\times\mathbb{R}^{(a,b)}$ and let $g_N$ be the product metric on $N$.
\begin{enumerate}
\item For generic $f$, $(N,g_N)$ is not locally homogeneous if $p\ge2$.
\item The manifold $(N,g_N)$ is not mixed Jordan IP.
\item Suppose that $a>0$ and that $b=0$. Then $(N,g_N)$ is neither timelike Jordan Osserman nor timelike Jordan Ivanov-Petrova.
Furthermore $(N,g_N)$ is spacelike Jordan Osserman and spacelike Jordan IP.
\item Suppose that $a=0$ and that $b>0$. Then $(N,g_N)$ is timelike Jordan Osserman and timelike Jordan Ivanov-Petrova.
Furthermore
$(N,g_N)$ is neither spacelike Jordan Osserman nor spacelike Jordan IP.
\item Suppose that $a>0$ and that $b>0$. Then $(N,g_N)$ is neither timelike Jordan Osserman nor timelike Jordan Ivanov-Petrova nor
spacelike Jordan Osserman nor spacelike Jordan Ivanov-Petrova.
\end{enumerate}
\end{theorem}

One also has \cite{GIZ02}:

\begin{theorem}\label{thm-3.x} Assume that $H$ is definite. Then $(N,g_N)$ is Jordan Osserman\begin{enumerate}
\item of types $(r,0)$ and $(p-r,q+b)$ if $a=0$ and if $0<r\le p$;
\item of types $(0,s)$ and $(p+a,q-s)$ if $b=0$ and if $0<s\le p$;
\item of types $(r,0)$ and $(p+a-r,q+b)$ if $a>0$ and if $a+2\le r\le p+a$;
\item of types $(0,s)$ and $(p+q,q+b-s)$ if $b>0$ and if $b+2\le s\le q+b$.
\end{enumerate} $(N,g_N)$ is not Jordan Osserman for other values of $(r,s)$.
\end{theorem}

\subsection{Manifolds of signature $(s,2s)$}\label{subs-3.4} The second family arises in signature $(2s,s)$ for $s\ge2$.
 Let $\vec u:=(u_1,...,u_s)$, $\vec t:=(t_1,...,t_s)$, and $\vec v:=(v_1,...,v_s)$ give
coordinates $(\vec u,\vec t,\vec v)$ on $\mathbb{R}^{3s}$ for $s\ge2$. Let 
$$F(\vec u):=f_1(u_1)+...+f_s(u_s)$$ 
be a smooth function
on an open subset $\mathcal{O}\subset\mathbb{R}^s$.
Define a
pseudo-Riemannian metric
$g_F$ of signature
$(2s,s)$ on $M:=\mathcal{O}\times\mathbb{R}^{2s}$ whose non-zero components are:
$$
\begin{array}{l}
g_F(\partial_i^u,\partial_i^u)=-2F(\vec u)-2\textstyle\sum_{1\le i\le s}u_it_i,\\
 g_F(\partial_i^u,\partial_i^v)=g_F(\partial_i^v,\partial_i^u)=1,\\
 g_F(\partial_i^t,\partial_i^t)=-1\,.
\end{array}$$

We also define the corresponding model spaces. Let 
$$\{U_1,...,U_s,T_1,...,T_s,V_1,...,V_s\}$$
 be a basis for $\mathbb{R}^{3s}$. Let $\mathcal{V}_{3s}:=(\mathbb{R}^{3s},g_V,A_V)$,
where the non-zero entries of the inner product $g_V$ and of the algebraic curvature tensor
$R_V$, up to the usual $\mathbb{Z}_2$ symmetries, are:
$$
\begin{array}{l}
g_V(U_i,V_i)=1,\quad
g_V(T_i,T_i)=-1,\quad\text{and}\\
R_V(U_i,U_j,U_j,T_i)=1\quad\text{for}\quad i\ne j\,.
\end{array}
$$
One then has \cite{GN03}:
\begin{theorem}\label{thm-3.4} The manifolds $(M,g_F)$ are curvature homogeneous with model $\mathcal{V}_{3s}$. They are not
locally homogeneous for generic $F$.
\end{theorem}

It is again immediate from an examination of the model space that
\begin{eqnarray*}
&&R(\cdot,\cdot):\operatorname{Span}\{U_i\}\rightarrow\operatorname{Span}\{T_i,V_i\},\\
&&R(\cdot,\cdot):\operatorname{Span}\{T_i\}\rightarrow\operatorname{Span}\{V_i\},\\
&&R(\cdot,\cdot):\operatorname{Span}\{V_i\}\rightarrow\{0\}\,.
\end{eqnarray*}
Consequently, the curvature operator is $3$-nilpotent. As the Ricci operator is nilpotent, the manifold $(M,g_F)$ is
Ricci flat and Einstein. One also has \cite{GN03}:

\begin{theorem}\label{thm-3.6} We have that $(M,g_F)$ is: \begin{enumerate}
\item spacelike Jordan Osserman but not timelike Jordan Osserman;
\item $k$-spacelike higher order Jordan Osserman for $2\le k\le s$;
\item  $k$-timelike higher order Jordan Osserman if
and only if
$s+2\le k\le 2s$;
\item spacelike Jordan Ivanov-Petrova;
\item neither timelike nor mixed Jordan Ivanov-Petrova;
\item spacelike Jordan $k$-Stanilov for $2\le k\le s$;
\item timelike Jordan $k$-Stanilov if and only if $k=2s$.
\end{enumerate}
\end{theorem}

Since these manifolds are Ricci flat, $W=R$. Consequently, we have as well that
\begin{theorem}\label{thm-3.xx} We have that $(M,g_F)$ is: \begin{enumerate}
\item conformally spacelike Jordan Osserman but not conformally timelike Jordan Osserman;
\item conformally spacelike Jordan Ivanov-Petrova;
\item neither conformally timelike nor conformally mixed Jordan Ivanov-Petrova.
\end{enumerate}
\end{theorem}

We remark that the rank of the skew-symmetric curvature operator is $4$; these are the only known examples of spacelike Jordan
Ivanov-Petrova manifolds which have rank $4$.

\section*{Acknowledgments} Research of P. Gilkey partially supported by the
MPI (Leipzig). Research of S. Nik\v cevi\'c partially supported by the Dierks Von Zweck
Stiftung (Essen), DAAD (Germany) and MM 1646 (Srbija). The second and fourth authors thank the
Technical University of Berlin where some of the research reported here was conducted.

\end{document}